\input boxedeps.tex 
\SetepsfEPSFSpecial 
\HideDisplacementBoxes
\def\figin#1#2{\medbreak
$$
 {\BoxedEPSF{#1.eps scaled
#2}%
}%
$$
\medbreak\noindent}
\def\eqref#1{(\ref{#1})}
\def\capt#1#2{\centerline{\ninepoint Figure #1. #2}\medbreak}
\documentstyle{amlts}
\begin{document}
\annalsline{151}{2000}
\received{December 17, 1998}
\startingpage{309}
\def\bye{\end{document}}
 \font\tenrm=cmr10
\catcode`\@=11
\font\twelvemsb=msbm10 scaled 1100
\font\tenmsb=msbm10
\font\ninemsb=msbm10 scaled 800
\newfam\msbfam
\textfont\msbfam=\twelvemsb  \scriptfont\msbfam=\ninemsb
  \scriptscriptfont\msbfam=\ninemsb
\def\msb@{\hexnumber@\msbfam}
\def\Bbb{\relax\ifmmode\let\next\Bbb@\else
 \def\next{\errmessage{Use \string\Bbb\space only in math
mode}}\fi\next}
\def\Bbb@#1{{\Bbb@@{#1}}}
\def\Bbb@@#1{\fam\msbfam#1}
\catcode`\@=12

 \catcode`\@=11
\font\twelveeuf=eufm10 scaled 1100
\font\teneuf=eufm10
\font\nineeuf=eufm7 scaled 1100
\newfam\euffam
\textfont\euffam=\twelveeuf  \scriptfont\euffam=\teneuf
  \scriptscriptfont\euffam=\nineeuf
\def\euf@{\hexnumber@\euffam}
\def\frak{\relax\ifmmode\let\next\frak@\else
 \def\next{\errmessage{Use \string\frak\space only in math
mode}}\fi\next}
\def\frak@#1{{\frak@@{#1}}}
\def\frak@@#1{\fam\euffam#1}
\catcode`\@=12


\newcommand{\C}{\Bbb C}
\newcommand{\R}{\Bbb R}
\newcommand{\Z}{\Bbb Z}
\newcommand{\N}{\Bbb N}
\newcommand{\Nn}{\{1,\ldots,n\}}
\newcommand{\ind}{{\rm ind}}
\newcommand{\lk}{{\rm lk}}
\newcommand{\const}{{\rm const}}

\newcommand{\tor}{S^1\times S^1}
\newcommand{\cp}{\C P^2}
\newcommand{\st}{{\rm st}}
\newcommand{\degree}{{\rm deg}}
\newcommand{\conj}{{\rm conj}}
\newcommand{\dd}{\partial}
\newcommand{\la}{\langle}
\newcommand{\ra}{\rangle}
\newcommand{\ignore}[1]{\relax}
 

\title{Real algebraic curves, the moment map\\ and amoebas}

 \shorttitle{Real algebraic curves, the moment map and amoebas}  

 \acknowledgements{Research was partially supported by the NSF (DMS\#9801726).}
\author{G. Mikhalkin}
\institutions{Harvard University,
Cambridge, MA\\
{\eightpoint {\it E-mail address\/}: mihalkin@math.harvard.edu}}

\bigbreak \centerline{\bf Abstract}
\bigbreak
In this paper we prove the topological uniqueness of
maximal arrangements of a real plane algebraic curve
with respect to three lines. More generally, we
prove the  topological uniqueness of  a maximally arranged
algebraic curve on a real toric surface.
We use the moment map as a tool for studying the
topology of real algebraic curves
and their complexifications.

\section{Introduction and   statement of results}

\demo{{\rm 1.1.\enspace M-}\/curves in the plane}
An algebraic curve $\R\bar{A}\subset\R P^2$ is the zero
set of a polynomial $p$ of degree $d$. (We reserve the notation
$\R A$ for the corresponding curve in $(\R-0)^2$.)
Suppose that $\R\bar{A}$ is nonsingular. Then it is homeomorphic
to a disjoint union of circles. 
If $d$ is even then each component of $\R\bar{A}$ bounds
a disk in $\R P^2$; such a component is called an {\it oval}.
If $d$ is odd then all the components but one are ovals and
the remaining one is a one-sided circle in $\R P^2$.
Harnack's inequality \cite{Ha} states that the number of components
is not greater than $\frac{(d-1)(d-2)}{2}+1$.
If it is equal to $\frac{(d-1)(d-2)}{2}+1$ then $\R\bar{A}$ is called an M-\/{\it curve}.
Let $l_1,\dots,l_n$ be lines in general position in $\R P^2$
(i.e. no three lines pass through the same point).
\enddemo

\demo{{D}efinition {\rm 1}}
We say that $\R\bar{A}$ is in {\it maximal position} with
respect to a collection of $n$ lines $l_1,\dots,l_n$
if $\R\bar{A}$ is an M-curve and there exist $n$ disjoint
arcs $a_1,\dots,a_n\subset\R\bar{A}$ such that $a_j$ intersect $l_j$
in $d$ points and all the arcs belong to the same component
of $\R\bar{A}$; see Figure~1.
\enddemo

\demo{{R}emark {\rm 1}}
The arcs from Definition~1 were called {\it bases of rank $1$}
in \cite{B}.
Brusotti used them to generalize the constructions of Harnack
\cite{Ha} and Hilbert \cite{Hi} and produce a larger variety 
of M-curves for all $d$.
A starting curve for Brusotti's construction is an M-curve
with 2 disjoint bases (possibly of higher

\figin{maxpos}{400}
\capt{1}{Maximal position with respect to lines.}

\noindent  rank).
The constructions of Harnack and Hilbert are special
cases where one takes a line or an ellipse (respectively)
for the starting curve. Note that both the line and
the ellipse are M-curves with arbitrarily many disjoint bases
of rank~$1$.
\enddemo

We call the topological type of $(\R P^2;\R\bar{A},l_1\cup\dots\cup l_n)$,
where $\R\bar{A}$ is a curve of degree $d$ in maximal position with
respect to $l_1\cup\dots\cup l_n$, a {\it maximal topological type}.
What are maximal topological types for given $n$ and $d$?

For $n=0$ this problem is a part of Hilbert's 16th problem
\cite{H}. It is an open question.
There are powerful theorems known which show that
the topological type of an M-curve is very restricted (\cite{P},
\cite{A}, \cite{R} et al.).
But on the other hand there is a large variety of
M-curves in different topological types.
The complete answer is known only for $d\le 7$; see \cite{V}.

For $n=1$ this asks for topological classification of
affine M-curves. Indeed, the maximality condition from
Definition~1 for $n=1$ just states that
the affine curve $\R\bar{A}-l_1\subset\R^2=\R P^2-l_1$
has $\frac{(d-1)(d-2)}{2}+d$ components, the maximal
possible value for a curve of degree $d$.
The affine M-curves of degree 5 were classified in \cite{Po};
see Figure~2.
The question is open for $d>5$.
\bigbreak
\begin{center}
\BoxedEPSF{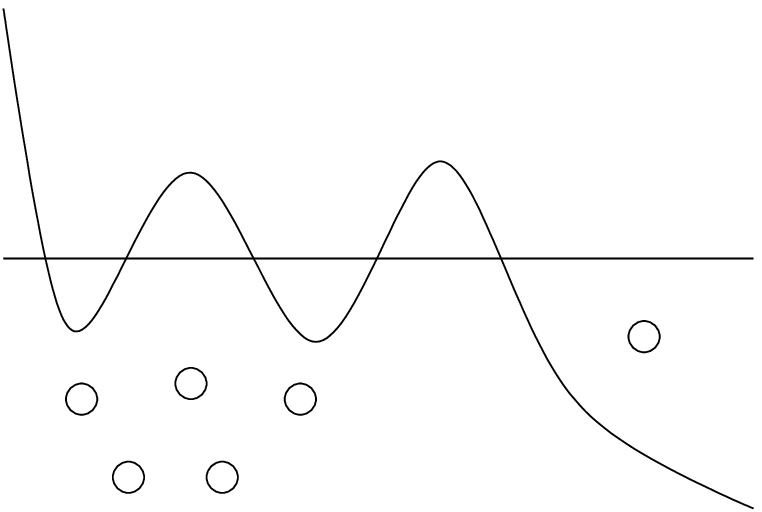 scaled 400}\hspace{.02in}\BoxedEPSF{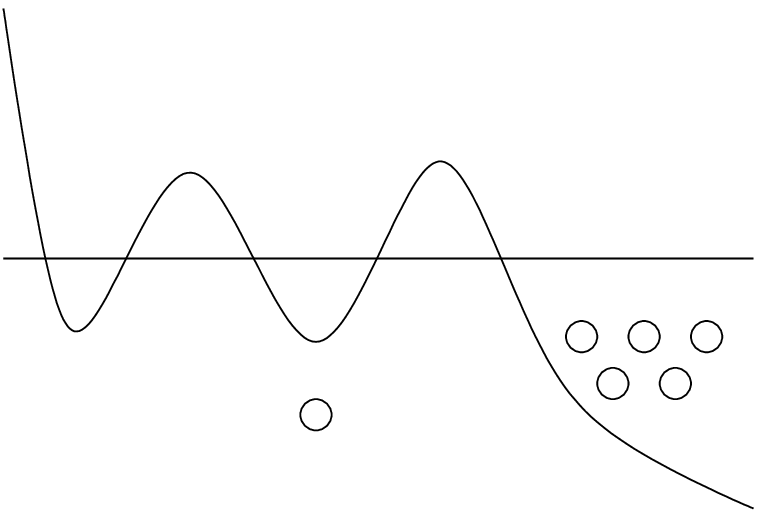 scaled 400}\hspace{.02in}
\BoxedEPSF{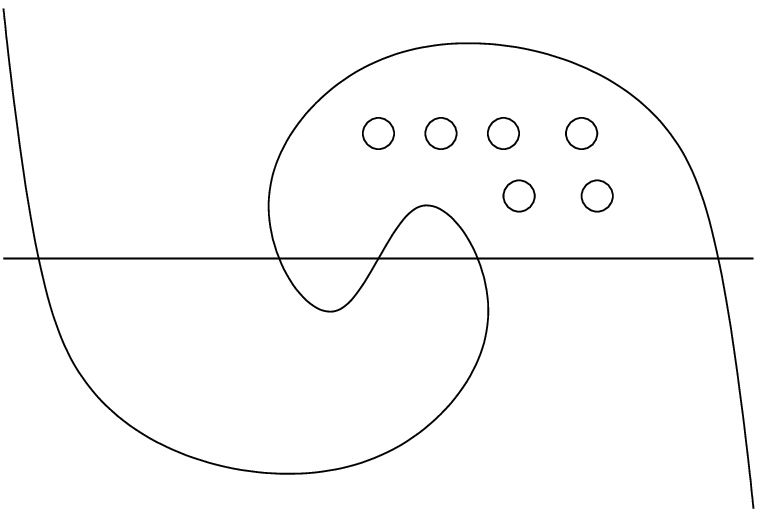 scaled 400}\hspace{0.2in}\BoxedEPSF{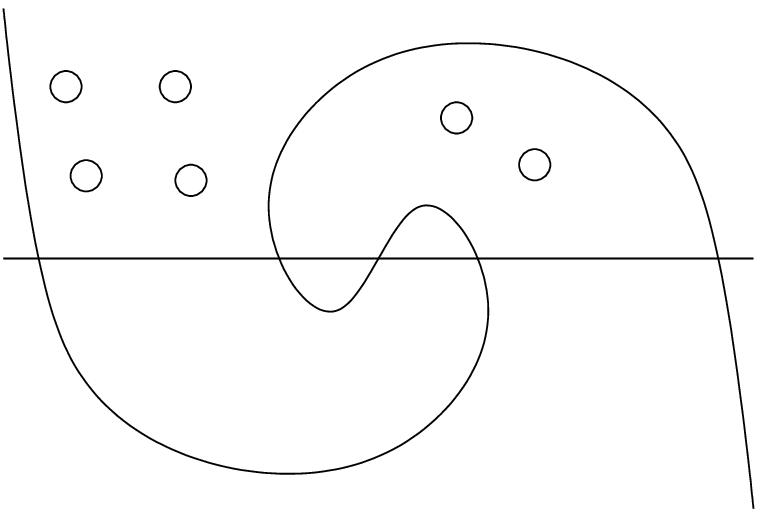 scaled 400}
\end{center}
\bigbreak
\capt{2}{Maximal topological types for $n=1$, $d=5$.}

For $n=2$ there are several constructions of M-curves
of the same degree $d$ but of different topological types
found by Brusotti \cite{B} with $d\ge 4$.
Figure~3 pictures all the types for $d=4$.
For $d>4$ the question is open.  

The following two theorems answer this question for all $d$ if $n\ge 3$.
For $n=3$ we can assume that $l_1,l_2,l_3$ are coordinate lines in $\R P^2$
(the $x$-axis, the $y$-axis and the infinite line).

\centerline{\BoxedEPSF{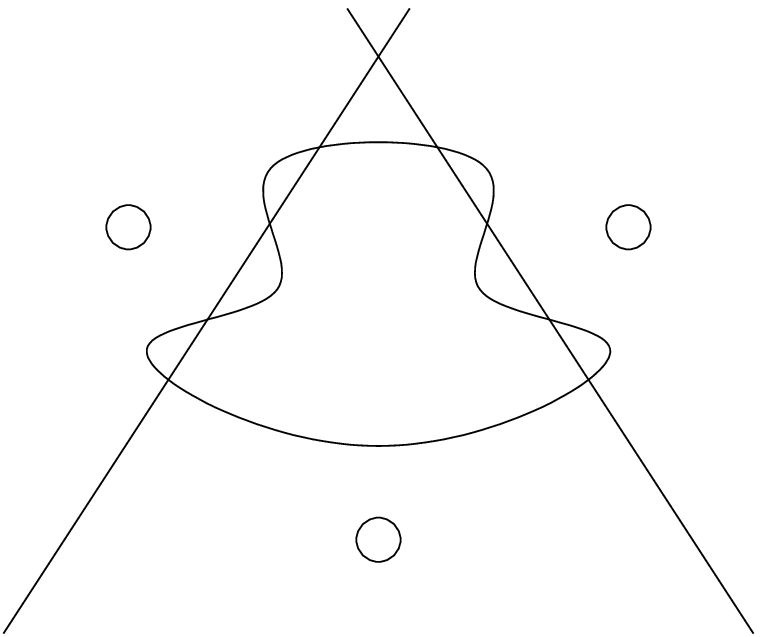 scaled400}\hspace{.5in}
\BoxedEPSF{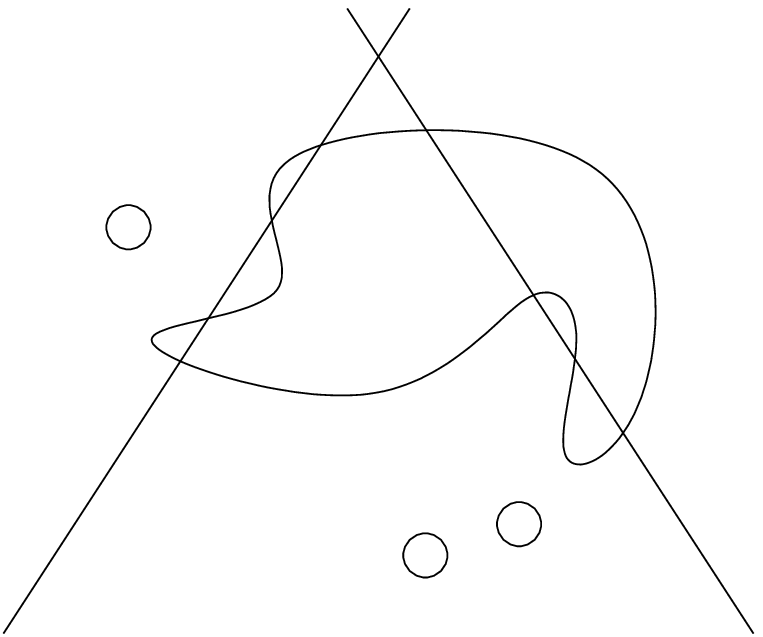 scaled400}\hspace{.5in}
\BoxedEPSF{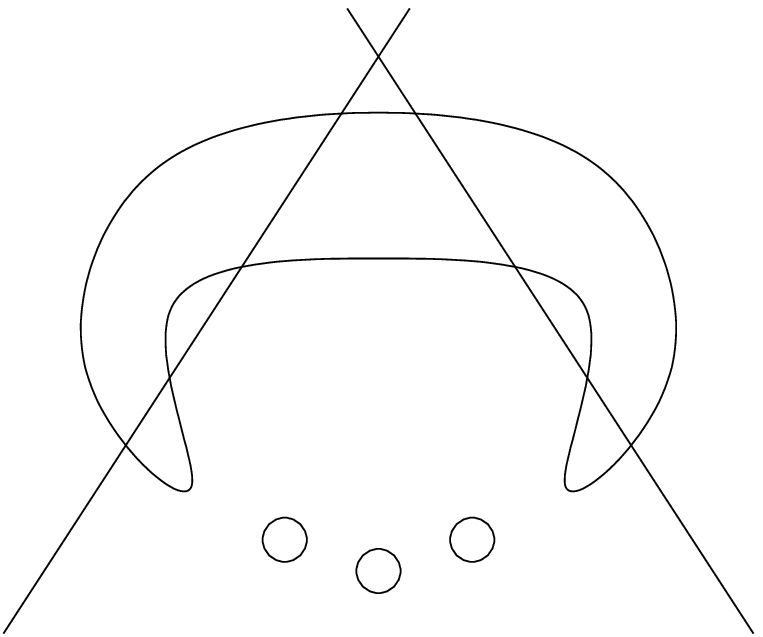 scaled400}}
\bigbreak
\capt{3}{Maximal topological types for $n=2$, $d=4$.}

\specialnumber{1} 
\proclaim{Theorem} \hskip-3pt
The maximal topological type is unique for $n=3$ and any~$d${\rm .}
If $\R\bar{A}\subset\R P^2$ is a curve
of degree $d$ in maximal position with respect to the coordinate
axes then the topological type of $(\R P^2;\R\bar{A},l_1\cup l_2\cup l_3)$
depends only on $d$ and is pictured in Figure~{\rm 4} for even $d$
or Figure~{\rm 5} for odd $d${\rm .}
\endproclaim

The M-curves pictured in Figure~{\rm 4} and Figure~{\rm 5} were
constructed by Harnack \cite{Ha}.
In fact, they were the first examples of
M-curves in $\R P^2$ for arbitrary $d$.

\specialnumber{2} \proclaim{Theorem}
There are no maximal topological types for $n>3$ and $d\ge 3${\rm .}
\endproclaim

In particular, this takes care of the case $n>4$ where
the choice of the lines
$l_1,\dots,l_n\subset\R P^2$ is not unique so that the answer might
presumably depend on it.

\medbreak
\centerline{\BoxedEPSF{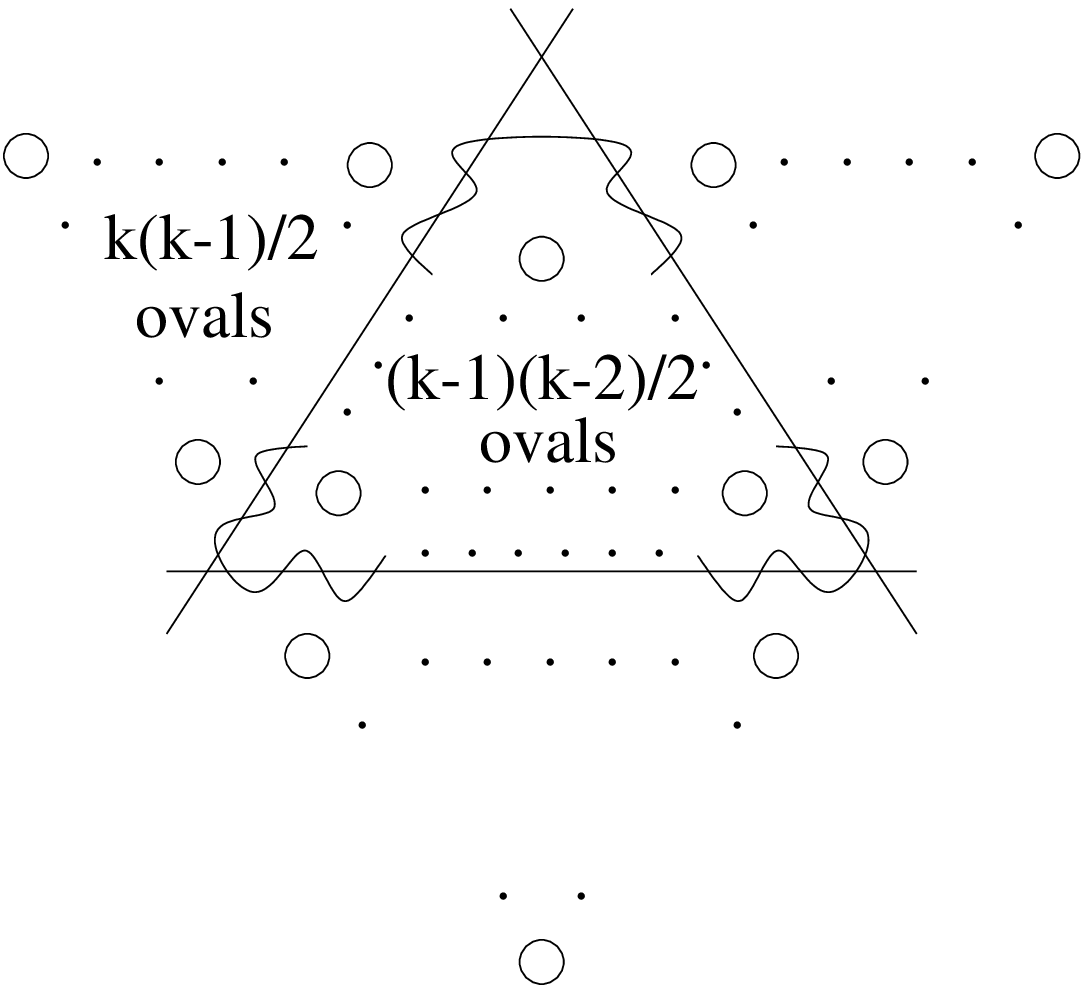 scaled 350}\hspace{.3in}
 \BoxedEPSF{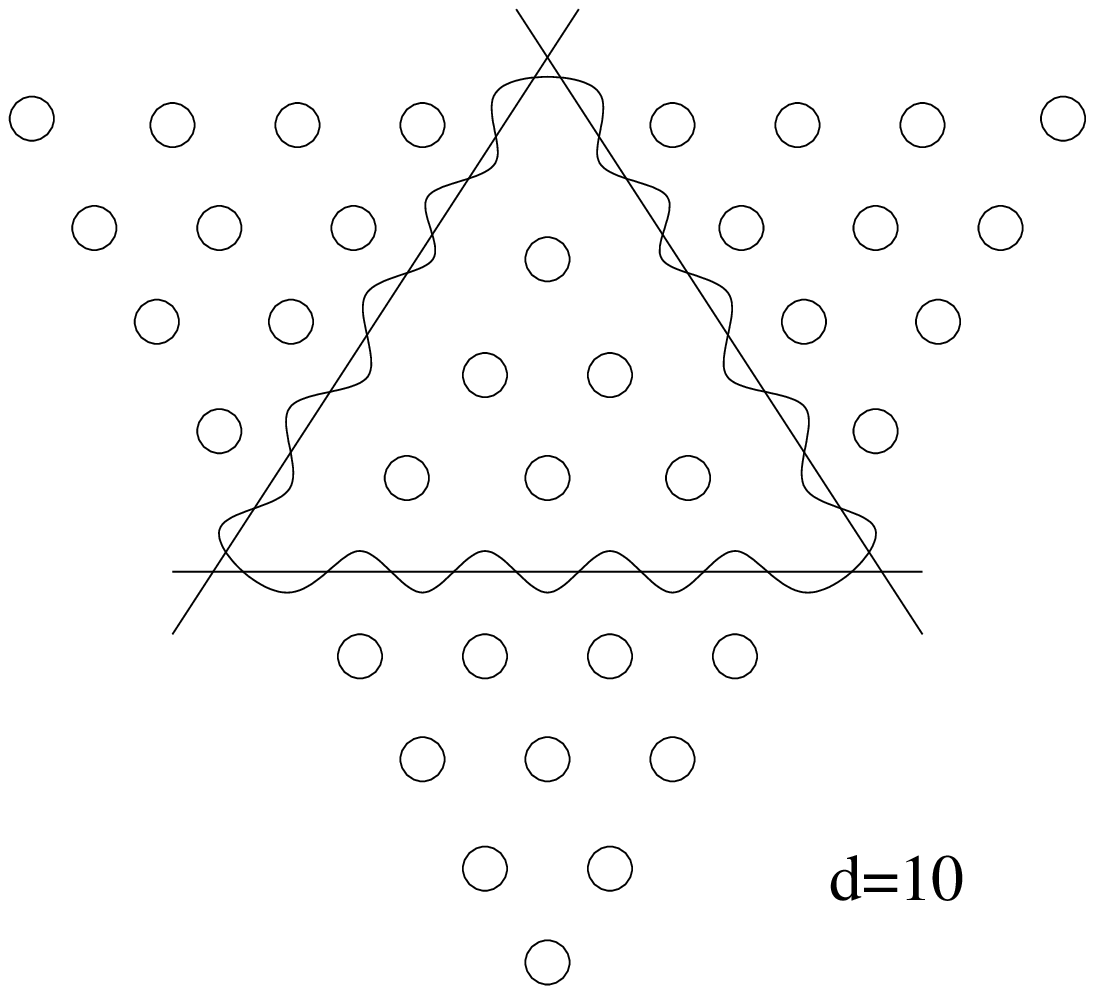 scaled 350}
}
\medbreak
\capt{4}{The maximal type
for $n=3$, $d=2k$; e.g.\ $d=10$.}
\smallbreak
\centerline{ \BoxedEPSF{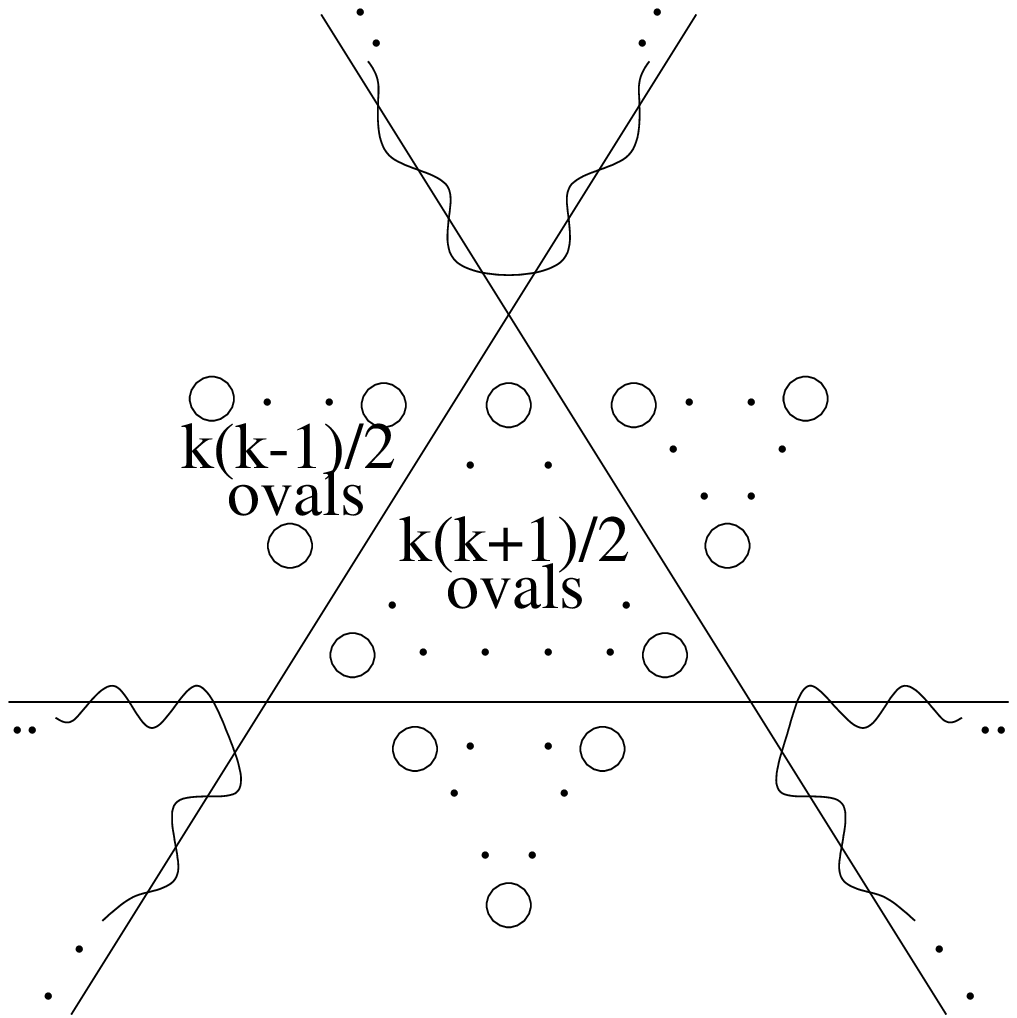 scaled 330}\hspace{.3in}
\BoxedEPSF{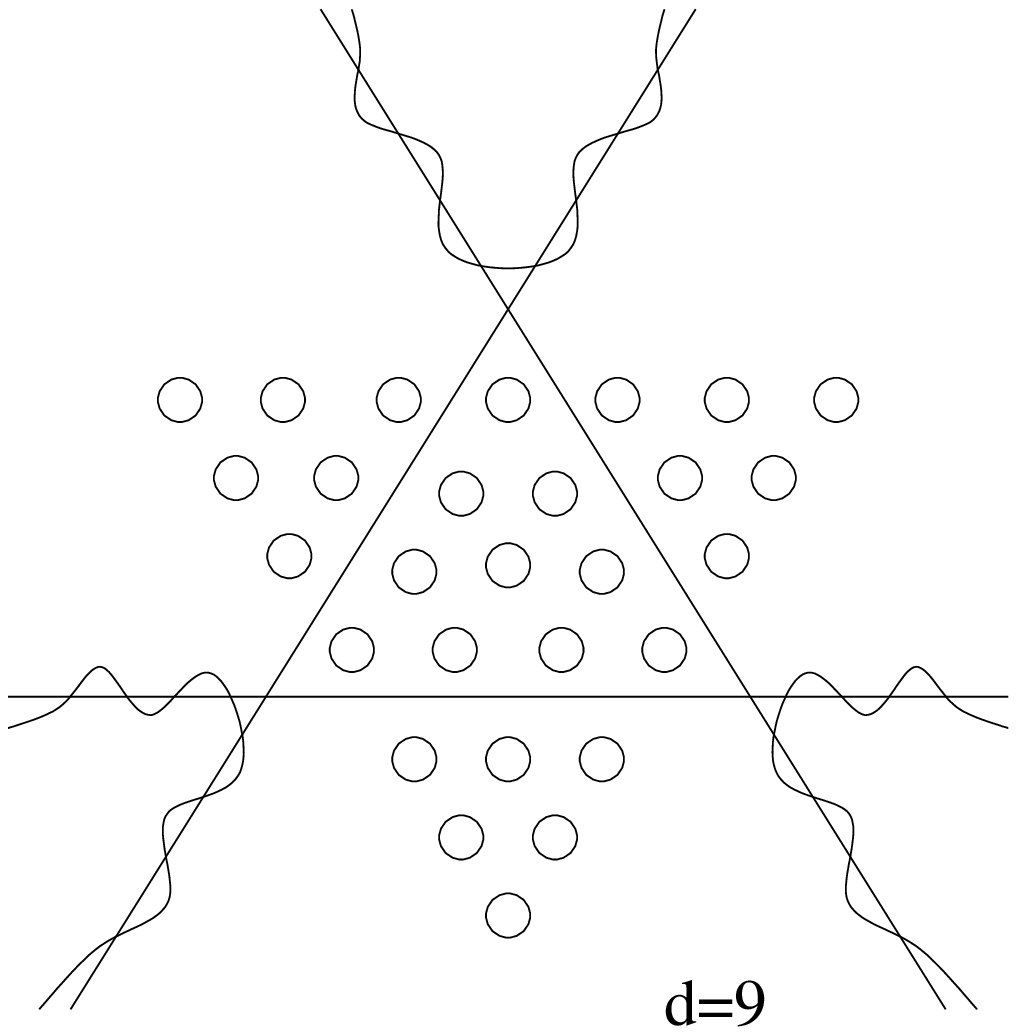 scaled 330}}
\capt{5}{The maximal type for $n=3$, $d=2k+1$; e.g.\ $d=9$.}
 
\demo{{R}emark {\rm 2}}
For every $d$ there exists an algebraic curve of degree $d$ in
maximal position with respect to three lines and invariant
with respect to an $S_3$ group of symmetries 
of these three lines. Such a curve can be constructed by patchworking;
see the appendix.
Note that the curves in Figure~{\rm 4} intersect each of the lines
in an infinite point (the topological arrangements pictured
there also admit an $S_3$-symmetry).
\enddemo
\demo{{R}emark {\rm 3}}
Let $d=3$. An arc of a real cubic curve is a base of rank~$1$
if and only if it contains an inflection point.
Thus, Theorem 2 may be viewed as a generalization of the
fact that a plane cubic curve cannot have more than three real inflection points.
\enddemo
\demo{{R}emark {\rm 4}}
The condition that $\R\bar{A}$ is an M-curve is essential for Theorem~2.
There exists a quartic with three ovals and with four bases
of rank~$1$ on the same oval.
The condition that the bases of rank~$1$ be on the same component is
also essential.
There exists an M-quartic with four bases of rank~$1$ on different ovals.
Such quartics may be obtained by perturbing a union of
four lines; see Figure~6.
The first quartic consists of three ovals and
two of them intersect the infinite line.
The second quartic consists of four ovals, none
intersects the infinite line.
\enddemo

\centerline{\BoxedEPSF{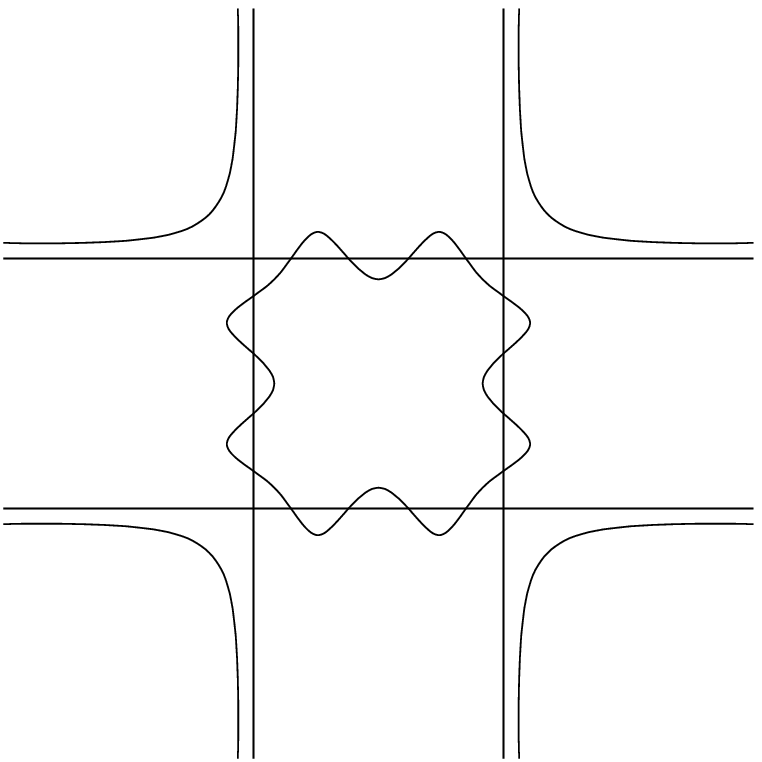 scaled 400}\hspace{1in}
\BoxedEPSF{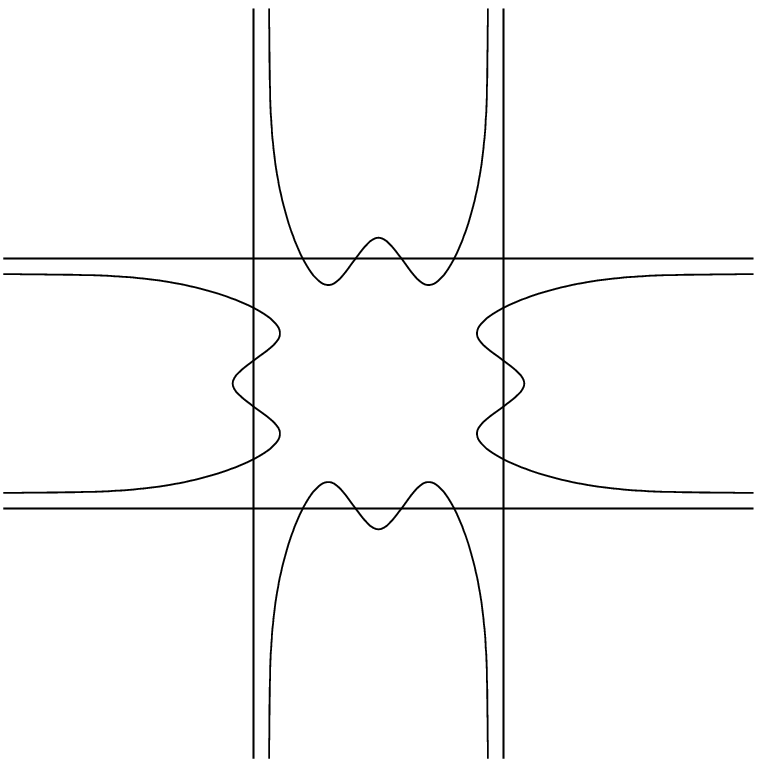 scaled 400}}
\medbreak
\capt{6}{Quartics with  four bases.}

\demo{{\rm 1.2. \enspace M-}\/curves in toric surfaces}
A convex polygon $\Delta\subset\R^2$
with vertices in $\Z^2$ defines a compactification
of $(\C-0)^2$ to an algebraic surface $T$ called
the toric surface; see e.g.\ \cite{GKZ}.
The action of the torus $S^1\times S^1$ on $(\C-0)^2$ by
$(\alpha,\beta)\times(x,y)=(\alpha x,\beta y)$,
where we view each $S^1$ as the unit circle in $\C$,
extends to the action of $S^1\times S^1$ on $T$.
The {moment map} $\mu_T:T\to\Delta$ of this action
exhibits $T$ as a singular Lagrangian fibration over $\Delta$.
The projection $\mu_T$ takes the quotient by the
action of the group $S^1\times S^1$.
The closure $\R T$ of $(\R-0)^2$ in $T$ is {the real toric
surface} and $T$ is its complexification.
Over the interior of $\Delta$ the fiber of $\mu_T$
is the full torus $S^1\times S^1$.
Over the sides of $\Delta$ the fiber is the quotient
of the torus by the circle subgroup corresponding to the
(rational) slope of the side. Over the vertices of $\Delta$
the fiber is one point.
The inverse images of the sides of $\Delta$
are rational holomorphic curves in $T$;
we call them {\it the axes} of $T$.

\demo{{E}xample}
Let $\Delta$ be the triangle
with vertices $(d,0)$, $(0,d)$ and $(0,0)$.
For any $d$, the corresponding toric surfaces
are the real projective plane $\R P^2$ and its
complexification $\C P^2$ (different values of $d$ correspond
to different multiples of the K\"ahler form).
The inverse images of the sides $[(0,0),(d,0)]$,
$[(0,0),(0,d)]$ and $[(d,0),(0,d)]$ are the $x$-axis,
the $y$-axis and the infinite line respectively.
\enddemo

Let $\R\bar{A}\subset\R T$ be an algebraic curve and $\bar{A}\subset T$
be its complexification.
It is given by a real polynomial $p$ in two variables.
Recall that the Newton polygon $\Delta$
of $p=\sum a_{j,k}x^jy^k$ is the convex hull
of $\{(j,k)\ |\ a_{(j,k)}\neq 0\}$ in $\R^2$.
If $\bar{A}$ does not pass through the intersection of
the axes of $T$ but intersects every axis,
then $T$ is necessarily the toric
surface corresponding to the Newton polygon $\Delta$.
If $\bar{A}$ is nonsingular then the genus $g$ of $\bar{A}$ is equal to
the number of lattice points in the interior of $\Delta$;
see \cite{Kh}. We call $\R\bar{A}$ an M-\/{\it curve} if the
number of components of $\R\bar{A}$ is equal to $g+1$;
by Harnack's inequality this is the largest possible number.
Let $l_1,\dots,l_n$ be the axes of $T$ in the order corresponding
to the order of the sides of $\Delta$. Let $d_j$ be the integer length
of $l_j$ (one plus the number of integer points inside $l_j$).
Note that $d_j$ is the degree of the restriction of $p$ to
$l_j$ and, therefore, the number of intersection points
of $l_j$ and $\R\bar{A}$ is no more than $d_j$.
\demo{{D}efinition {\rm 2}}
We say that $\R\bar{A}$ is in {\it maximal position} in $\R T$
if $\R\bar{A}$ is an M-curve and there exist $n$ disjoint
arcs $c_1,\dots,c_n\subset\R\bar{A}$ such that the  $c_j$ intersect $l_j$
in $d_j$ points and all the arcs belong to the same component
of $\R\bar{A}$.
We say that $\R\bar{A}$ is in {\it cyclically maximal position} in $\R T$
if, in addition, the order of arcs $c_j$ agrees with the cyclic order
on the component of $\R\bar{A}$.
\enddemo
\demo{{R}emark {\rm 5}}
Obviously, if $n\le 3$ then a curve in maximal position in $\R T$
is automatically in cyclically maximal position in $\R T$.
The same is true if $n>3$ and all $d_j$ are even.
See Example~1 in Section \ref{realamoebas} for
a curve in maximal but not in cyclically maximal position
in $\R T$.
\enddemo
The following theorem is a generalization of Theorem 1.
\specialnumber{3} 
\proclaim{Theorem}
If $\R\bar{A}$ is in cyclically maximal position in $\R T$ then
the topological type of $(\R T;\R\bar{A},l_1\cup\dots\cup l_n)$
depends only on $\Delta${\rm .}
\endproclaim

The maximal topological type for $\Delta$ can be reconstructed from
Figure~12 and Lemma \ref{keylem}.
\demo{{R}emark {\rm 6}} The maximality definition and Theorem 3
can be restated in terms of curves in $(\R-0)^2$.
\enddemo

Curves in maximal position in $\R T$ exist for any $\Delta$
and the corresponding toric surface $T$; see Corollary A4
of the appendix.
\demo{{R}emark {\rm 7}}
For some shapes of $\Delta$ the surface $T$ has singularities
at the intersection of axes but (by the maximal position assumption)
$\R\bar{A}$ does not pass through them.
\enddemo

\demo{{R}emark {\rm 8}}
Note that different polygons $\Delta$ with parallel sides
correspond to different homology classes of $\bar{A}$ in the same $T$.
\enddemo

\section{Amoebas}

In this section we introduce the tools for proving the
main results. They are based on the following notion of
{\it amoeba} defined by Gelfand et al. in \cite{GKZ}.
\demo{{D}efinition {\rm 3}}
Let $p:\C^m\to\C$ be a polynomial.
Its amoeba is the image $\mu(A)\subset\R^m$,
where $\mu:(\C-0)^m\to\R^m$ is given by the formula
$$\mu(z_1,\dots,z_m)=(\log|z_1|,\dots,\log|z_m|)$$
and $A\subset(\C-0)^m$ is the zero locus of $p$
in the complement of the coordinate hyperplanes.
\enddemo
Let $\Delta$ be the Newton polytope
of $p$, $T\supset(\C-0)^m$ be the toric variety associated to $\Delta$
and $\mu_T:T\to\Delta$ be the moment map of $T$;
see e.g.\ \cite{GKZ}.
Note that $\mu:(\C-0)^m\to\R^m$ is just a reparametrization
of $\mu_T|_{(\C-0)^m}:(\C-0)^m\to{\rm Int}\Delta\subset\R^m$,
where the interior of $\Delta$ gets mapped to the whole $\R^m$.
Let $\bar{A}$ be the closure in $T$ of $A\subset(\C-0)^m\subset T$.
The image $\mu_T(\bar{A})\subset\Delta$ is called the compactified amoeba
\cite{GKZ}.
\medbreak
\centerline{\BoxedEPSF{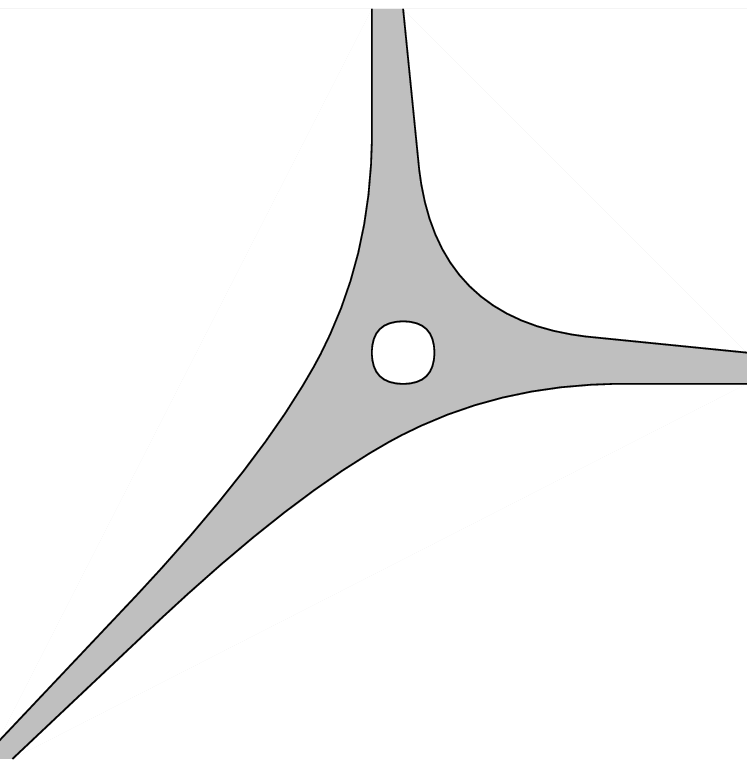 scaled 400}\hspace{1in}
\BoxedEPSF{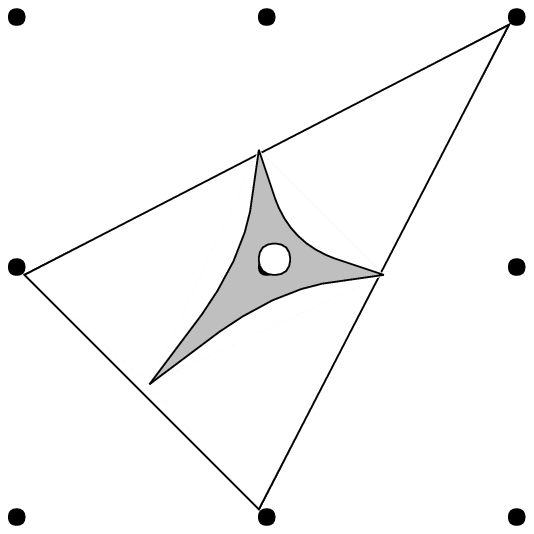 scaled 400}}
\medbreak
\capt{7}{Amoeba $\mu(A)$ and compactified amoeba $\mu_T(\bar{A})$.}
\smallbreak

Consider the region $\R^m-\mu(A)$.
It consists of bounded and unbounded components.
It was observed in \cite{GKZ} that each component
of $\R^m-\mu(A)$ is convex and that unbounded components
of $\R^m-\mu(A)$ correspond (inductively) to the complementary
regions of the amoebas of the intersection of $\bar{A}$ with
the toric subvarieties corresponding to the faces of $\Delta$.

\demo{{R}emark {\rm 9}} For $m=2$ these intersections are zero-dimensional,
so that  the description of the unbounded components is easy.
Indeed, the {\it tentacles} (i.e. the ends of $\mu(A)$)
correspond to the intersection of $\bar{A}$ with the
axes $l_j$. The number of such points counted with multiplicities
is $d_j$, since $d_j$ is the degree of the restriction of $p$
to $l_j$. Thus, if $\bar{A}$ intersects $l_j$ transversely and
no two intersection points are on the same fiber of $\bar{\mu}$,
the unbounded components of $\R^2-\mu(A)$ are in one-to-one correspondence
with the integer points on the boundary of $\Delta$,
i.e. with the set $\dd\Delta\cap\Z^m$.
\enddemo

The number of components of $\R^m-\mu(A)$ depends
on the coefficients of $p$ and not just on $\Delta$.
However, Forsberg et al. \cite{FPT} obtained the following
upper bound for this number in terms of $\Delta$.

\proclaimtitle{\cite{FPT}}
\specialnumber{4} \proclaim{Theorem}
\label{comp}
There is a natural injective map $\ind$ from the set
of components of $\R^m-\mu(A)$ to $\Delta\cap\Z^n${\rm .}
\endproclaim

There do exist amoebas with a bijective map $\ind$; see the
appendix for the construction.
In the proof of our main theorems we use a topological
interpretation of $\ind$.

\section{Amoebas from a topological point of view}

In this section we give a topological proof
of Theorem \ref{comp}.
Also we define the {\it logarithmic Gauss map} $\gamma$ for $\bar{A}$
and compute its degree.

\demo{{\rm 3.1.}\enspace Proof of Theorem {\rm \ref{comp}}}
If $x\in\R^m-\mu(A)$ then $\mu^{-1}(x)$ is an m-dimensional torus
in $(\C-0)^m$ not intersecting the hypersurface $A$.
Therefore the linking number with the closure of
$A$ in $\C^m$ produces a well-defined linear function 
${\rm lk}:H_1(\mu^{-1}(A))\to\Z$.
But $H_1(\mu^{-1}(A))=\Z^m$, where the identification is given
by coordinates in $\C^m$. Therefore, $\lk$ is given by $m$ integer
numbers $(i_1,\dots,i_m)$.
Define $\ind:\R^m-\mu(A)\to\Delta\cap\Z^m$ by $\ind(x)=(i_1,\dots,i_m)$.
Clearly, $\ind$ is locally constant and therefore defines a map
on the set components of $\R^m-\mu(A)$. We need to prove that this
map is injective and lands on $\Delta\cap\Z^m$.

Denote by $\pi:\Z^m\to\Z$ the projection
defined by $$\pi_(j_1,\dots,j_m)=k_1j_1+\dots
+k_mj_m.
$$
 To prove that $\ind(x)\in\Delta$ it suffices to
prove that $\pi(\ind(x))\in\pi(\Delta)$
for any $(k_1,\dots,k_m)\in\Z^m$.

Let $C\subset\C^m$ be the curve given by the parametrization
$z_j=c_jt^{k_j}$, $0\neq c_j\in\C$.
For a generic choice of $c_j$, the pull-back of $p$ to $C$ is
a polynomial in one variable whose Newton polytope is
$\pi(\Delta)$.
But $\mu|_{C-(0,\dots,0)}$ is a circle fibration over the
line $x_j=k_jt+\log|c_j|$.
For $z\in C$ the circle $C\cap\mu^{-1}(\mu(z))$ represents
the homology class $(k_1,\dots,k_m)\in H_1(\mu^{-1}(\mu(z)))$.
The linking number in $\C^m$ of $C\cap\mu^{-1}(\mu(z))$ and $A$
is $d$, if $z\in C$ is sufficiently close to $(0,\dots,0)$, $D$,
if $z\in C$ is sufficiently close to infinity and anything in
between for other choices of $z$, where $[d,D]=\pi(\Delta)\subset\Z$.
This holds since we may use a part of $C$ as a membrane
to compute the linking number.
This implies that $\ind(x)\in\Delta$ since
the line $\mu(C)$ passes through the component of $\R^m-\mu(A)$
containing $x$ for a suitable choice of $c_j$.

To show the injectivity of ${\rm ind}$  choose a line $l\subset\R^n$
with a rational slope (i.e. given by $x_j=k_j+b_j$ for some
$k_j\in\Z$ and $b_j\in\R$)
which passes through any pair of components in $\R^m-\mu(A)$.
Let $x,y\in l$ be two points in different components of $\R^m-\mu(A)$.
Since $\mu^{-1}[x,y]\cap A\neq 0$ we may choose $c_j\in\C$ so that
$\mu(C)=l$ and $C\cap\mu^{-1}[x,y]\cap A\neq 0$ for $C$ parametrized
by $z_j=c_jt^{k_j}$. Therefore $\pi(\ind(x))\neq\pi(\ind(y))$ and
$\pi$ is injective.

Note that the same argument also proves
the convexity of components of $\R^m-\mu(A)$.
Convexity holds even locally as the following lemma
shows.
\specialnumber{1} \proclaim{Lemma}
\label{locconv}
Let $z\in A$ be a critical point of $\mu|_A${\rm ,}
$U\ni f(z)$ be a convex neighborhood of $f(z)$ in $\R^2$
and $V$ be the component of $(\mu_A)^{-1}(U)$ which contains $z${\rm .}
Then each component of $U-\mu(V)$ is convex{\rm .}
\endproclaim

\demo{Proof}
If not then there exists a straight closed interval
$I\subset U\subset\R^2$ with   both endpoints in the same
component of $U-\mu(V)$ which intersects $V$.
We may assume that $I$ has a rational slope and
that $I$ is transverse to $\mu|A$.
Since $U$ is contractible, $V\cap\mu^{-1}(I)$ is
null-homologous in $V$ and, therefore, null-homologous
in $(\C-0)^m$.
But since the slope of $I$ is rational there exists
a holomorphic annulus $Z$ which projects properly to $I$
and intersects $V$. This leads to a contradiction.
The intersection number of $Z$ and $V\cap\mu^{-1}(I)$
in $\mu^{-1}(I)$ is positive on one hand, since $V$ and $Z$
are holomorphic, and zero on the other hand, since
$V\cap\mu^{-1}(I)$ is null-homologous in $(\C-0)^m$ and,
therefore, null-homologous in $\mu^{-1}(I)$.
\enddemo

\demo{{\rm 3.2.}\enspace The logarithmic Gauss map {\rm (cf. \cite{Ka})}}
Suppose that $\bar{A}\subset T$ is nonsingular.
Let $\gamma:A\to\C P^{m-1}$ be the map defined
by $$\gamma(z_1,\dots,z_m)=[z_1\frac{\dd p}{\dd z_1}(z_1,\dots,z_m):
\dots:z_m\frac{\dd p}{\dd z_m}(z_1,\dots,z_m)].$$
The geometric description of $\gamma$ is the following.
Let $z\in A$ and $U\ni z$ be a small neighborhood of $z$ in $(\C-0)^m$.
Choose a branch of the holomorphic logarithm $\log_U:U\to\C^m$,
$(z_1,\dots,z_m)$ and apply the Gauss map $G$ to the image of $A$
(i.e. map each point $\log_U(A\cap U)$ to the tangent hyperplane
at that point).
The composition $G\circ\log_U$ does not depend on the
choice of the branch of the logarithm and gives $\gamma$.

Suppose that for any face of $\Delta$ the toric subvariety
corresponding to that face intersects $\bar{A}$ transversely.
In this case $\gamma$ extends to $\bar{A}$.
Let us extend $\gamma$ to $\bar{A}\cap
\mu^{-1}({\rm int}\Delta')$
for any facet $\Delta'\subset\Delta$.
Without   loss of generality we may assume that
the hyperplane corresponding to $\Delta'$ is $x_m=0$;
otherwise we change the coordinates in $(\C-0)^m$
by $Z_j=z_1^{b_{j1}}\dots z_m^{b_{jm}}$
for a  suitable integer $b_{jk}$.
Denote by $p_{\Delta'}=\sum\limits_{\Delta'}
a_{j_1,\dots,j_m}z_1^{j_1}\dots z_m^{j_m}$
the truncation of $p=\sum\limits_{\Delta}
a_{j_1,\dots,j_m}z_1^{j_1}\dots z_m^{j_m}$ to $\Delta'$.
Then $p_{\Delta'}$ is a polynomial in $(m-1)$ variables
$z_1,\dots,z_{m-1}$.
Define $$\bar\gamma(z_1,\dots,z_{m-1},0)=
[z_1\frac{\dd p_{\Delta'}}{\dd z_1}(z_1,\dots,z_{m-1}):
\dots:z_{m-1}\frac{\dd p_{\Delta'}}{\dd z_{m-1}}(z_1,\dots,z_{m-1}):0]$$
for $(z_1,\dots,z_{m-1},0)\in
\bar{A}\cap\mu^{-1}({\rm int}\Delta')$.
Inductively by codimension of the faces of $\Delta$,
$\gamma$ extends to a holomorphic map
$$\bar\gamma:\bar{A}\to\C P^{m-1}$$ between closed manifolds of
the same dimension.
\specialnumber{2} \proclaim{Lemma}
\label{gauss}
The degree of $\bar\gamma$ is $n!{\rm Vol}\Delta${\rm .}
\endproclaim
\demo{Proof}
The inverse image $\bar\gamma^{-1}([0:\dots:0:1])$ is the
zero set of $z_j\frac{\dd p}{\dd z_j}$,
$j=1,\dots,m-1$, and $p$.
The Newton polytope of each of these polynomials is $\Delta$.
By Kouchnirenko's theorem \cite{Kou} the number of points
in $\bar\gamma^{-1}([0:\dots:0:1])$
counted with multiplicities is $n!{\rm Vol}\Delta$.
\enddemo

\section{Real two-dimensional amoebas}
\label{realamoebas}

Suppose now that the coefficients of $p$ are real and $m=2$.
Then $\mu|_{{A}}\break :{A}\to\R^2$ is a map between smooth surfaces.
Its generic singularities in the class of smooth maps are
folds and cusps.
Denote by $F\subset A$ the locus of critical points
of $\mu|_A$, i.e. the points where $\mu_{{A}}$ is not
submersive. 

\specialnumber{3} \proclaim{Lemma}
\label{foldgauss}
${F}=\gamma^{-1}(\R P^1)${\rm .}
\endproclaim
\demo{Proof}
The holomorphic (multivalued) logarithm maps
the fiber tori\break $\mu^{-1}(x)$, $x\in\R^2$ to the purely imaginary planes
$\{{\rm Re}z_1=y_1,{\rm Re}z_2=y_2\}$
in $\C^2$.
Therefore, $z\in A$ is a critical point of $\mu_A$ if and only
if the logarithmic image of the tangent plane at $z$ contains
a purely imaginary vector. But this holds if and only if
$\gamma(z)\in\R P^1\subset\C P^1$.
\enddemo

Denote by $\R A\subset A$ the set of real zeroes of $p$ in $\R^2$
and by $\R\bar{A}\subset\bar{A}$ its compactification in $\R T$.

\specialnumber{4} \proclaim{{C}orollary}
$\R A\subset F${\rm .}
\endproclaim

Indeed, $\gamma(\R A)\subset\R P^1$ from the definition of $\gamma$.
However, $\R A$ does not always coincide with $F$;
e.g.\ $\R A$ may be empty while $F$, being the folds of a proper
degree 0 map $A\to\R^2$, is never empty.
If a point $z$ belongs to $F$ then its conjugate point $\bar{z}$
also belongs to $F$. Thus, imaginary folds are double folds and 
the number of points in the inverse image of $\mu$ jumps by four at those
folds.
\medbreak
\centerline{\BoxedEPSF{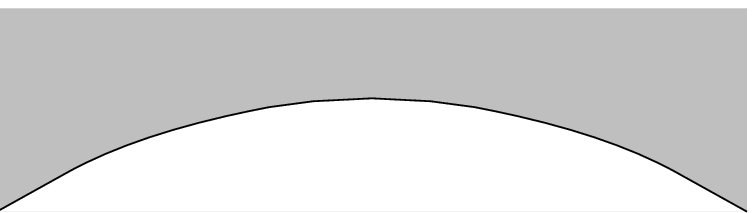 scaled 400}\hspace{1in}
\BoxedEPSF{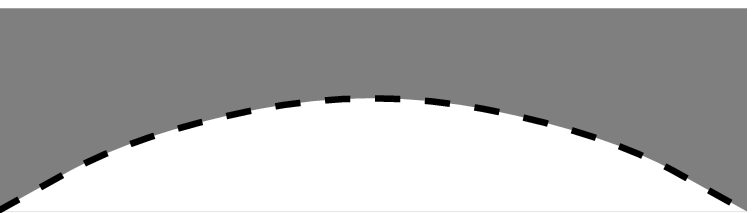 scaled 400}}
\medbreak
\capt{8}{Real folding and imaginary double folding.}

\demo{{R}emark {\rm 10}}
Besides folds and cusps, which are stable singularities
for smooth maps between surfaces there are two new stable
singularities for $\mu$ pictured in Figure~9.
These singularities persist under small real perturbations of $p$
but decompose into folds and cusps under small imaginary
perturbations.
Indeed, $\bar\gamma:\bar{A}\to\C P^1$ is a branched covering.
Branching points are the points of logarithmic inflection,
i.e. inflection after taking the holomorphic logarithm.
Generically, the branching is of multiplicity 2.
The new singularities correspond to the case when the branching
points are real. By Lemma \ref{foldgauss}, $F$ has double points
at the double branching points of $\gamma$. 
Only one of the branches at a double point of $F$ may be from $\R A$.
If the image of the other one under $\mu$ is not constant then
it corresponds to the imaginary double folds   pictured in
Figure~9 on the right. If it is constant then a whole circle
in $F$ is mapped to a point and this circle must have another real
point where the circle meet another branch of $\R A$. This is
a pinching singularity pictured in Figure~9 on the left.
It is a double point of $\mu(\R A)$ and both branches have an
inflection at this point.
The next two examples show that both cases appear.
\enddemo

\centerline{\BoxedEPSF{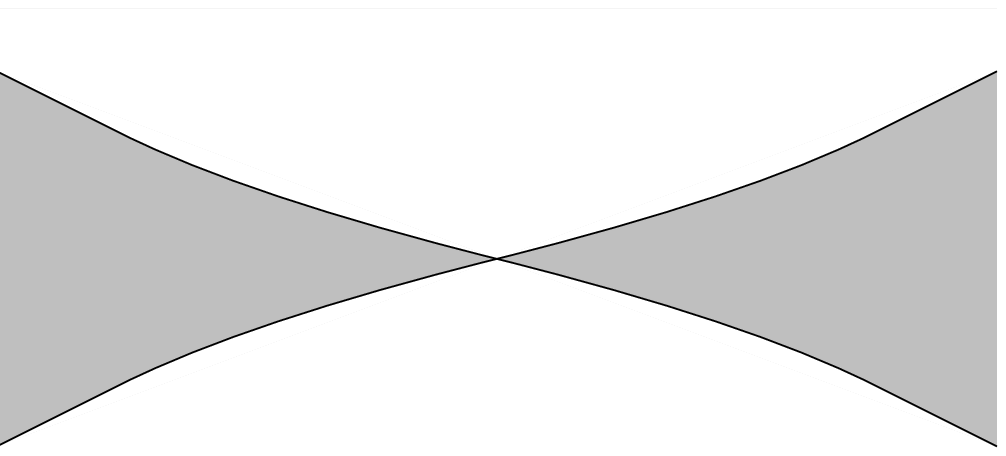 scaled 400}\hspace{1in}
\BoxedEPSF{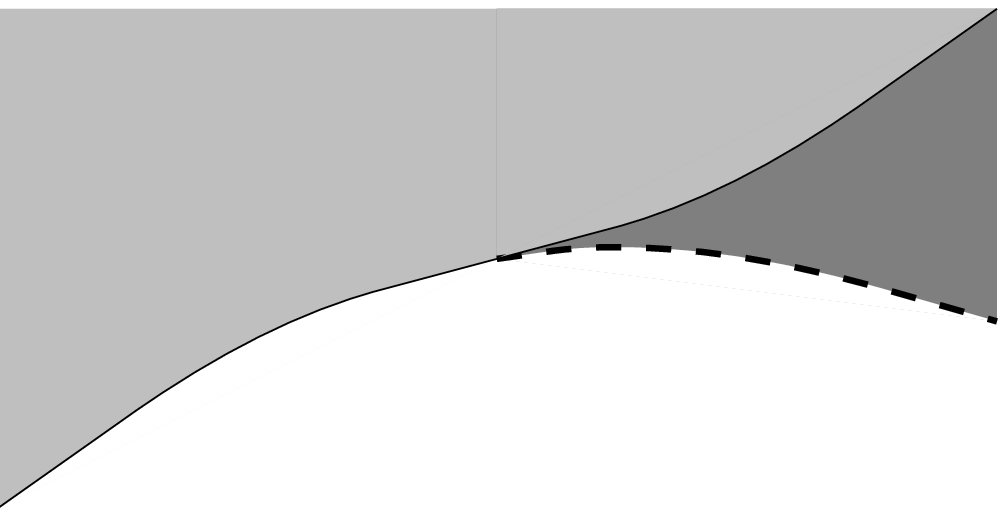 scaled 400}}
\medbreak
\capt{9}{A pinching and a junction of real
and imaginary folds.}

\demo{{E}xample {\rm 1}}
Let $p(x,y)=xy-x-y+a$, where $0<a<1$.
The corresponding hyperbola and its amoeba are pictured
in Figure~10.
\medbreak
\centerline{\BoxedEPSF{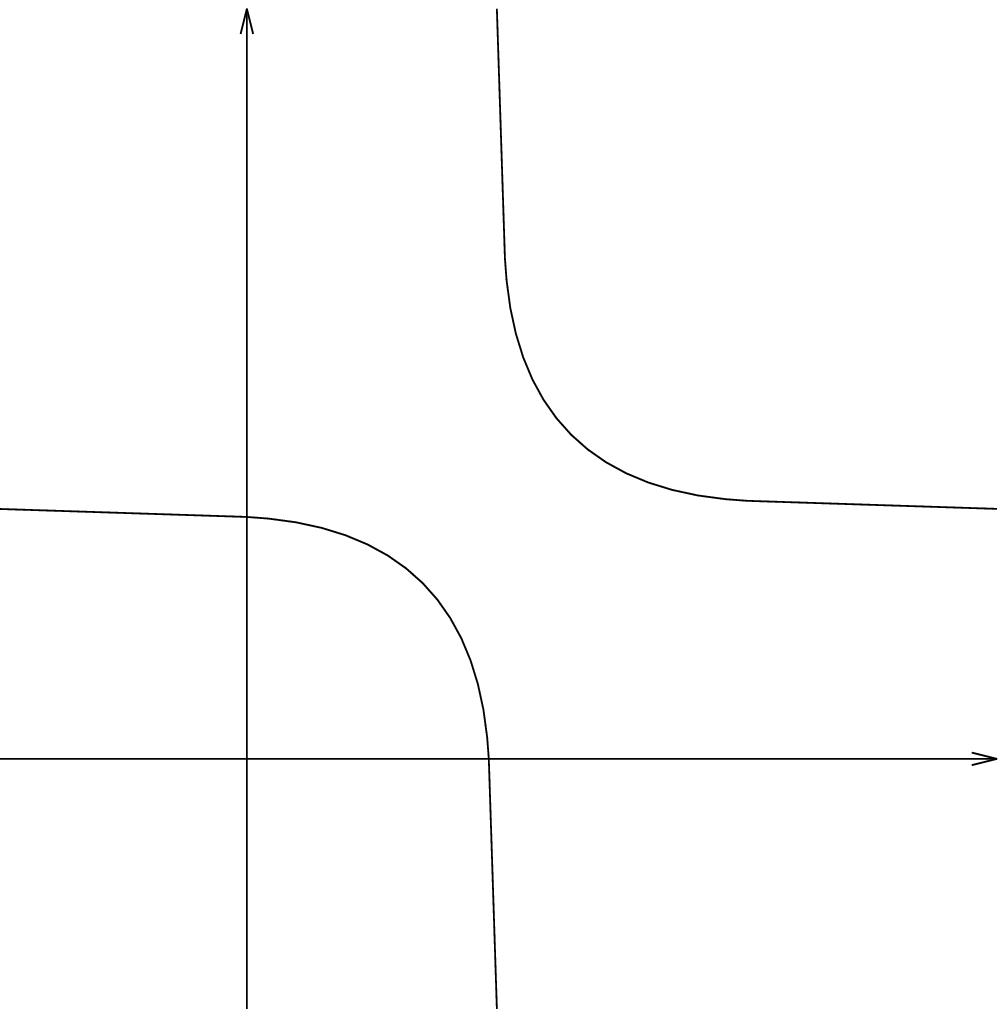 scaled 350}\hspace{1in}
\BoxedEPSF{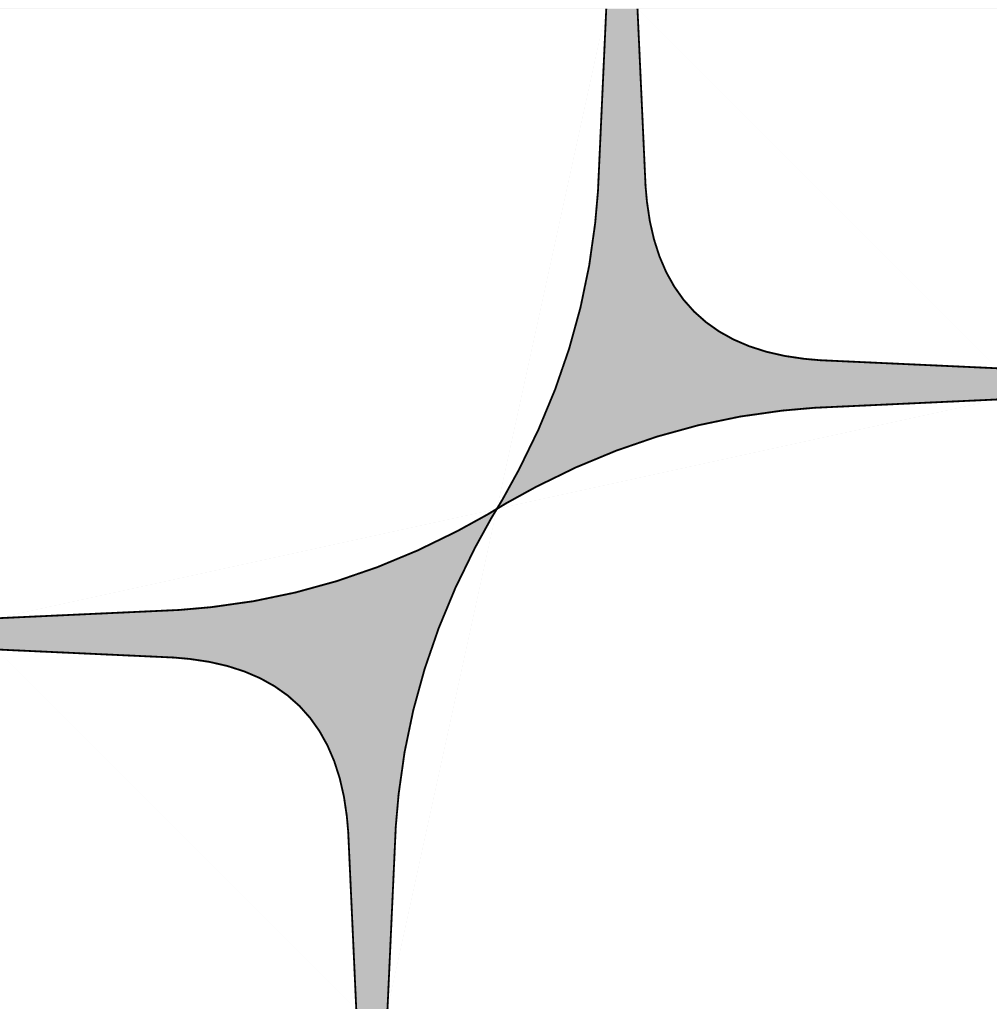 scaled 350}}
\vfill
\capt{10}{A hyperbola and its amoeba.}
\medbreak
Indeed, the curve $A$ is defined by $y-1=\frac{1-a}{x-1}$.
The image of the circle $|x|=c_x$ under $\frac{1-a}{x-1}$
is a circle which intersects the circle $|y|=c_y$ in two
points if $(c_x,c_y)\in\mu(A)-\mu(\R A)$, is tangent to it
if $(c_x,c_y)\in\mu(\R A)$, is disjoint to
it if $(c_x,c_y)\notin\mu(A)$ and coincides with it if
$c_x=c_y=\sqrt{a}$.
\enddemo

\demo{{E}xample {\rm 2}}
Let $p(x,y)=y-x^2+2x-a$, $a>1$.
The corresponding parabola and its amoeba are pictured
in Figure~11.
Indeed, $\R T$ in this case is a weighted projective space
with weights $(1,2,1)$ and $\R A$ intersects the $x$-axis
in a pair of conjugate imaginary points. The double imaginary
folding must merge to the real folding (by Lemma \ref{gauss}
the degree of $\gamma$ is 2).
\medbreak
\centerline{\BoxedEPSF{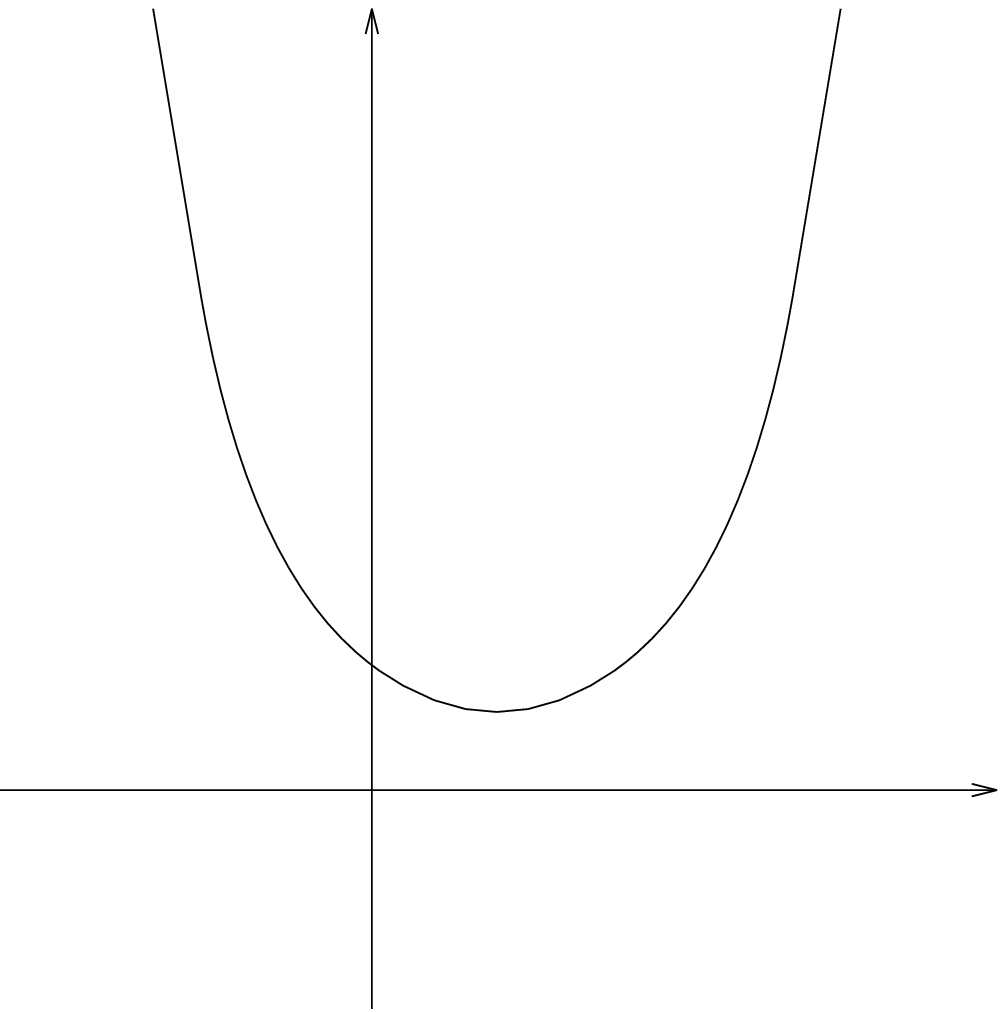 scaled 350}\hspace{1in}
\BoxedEPSF{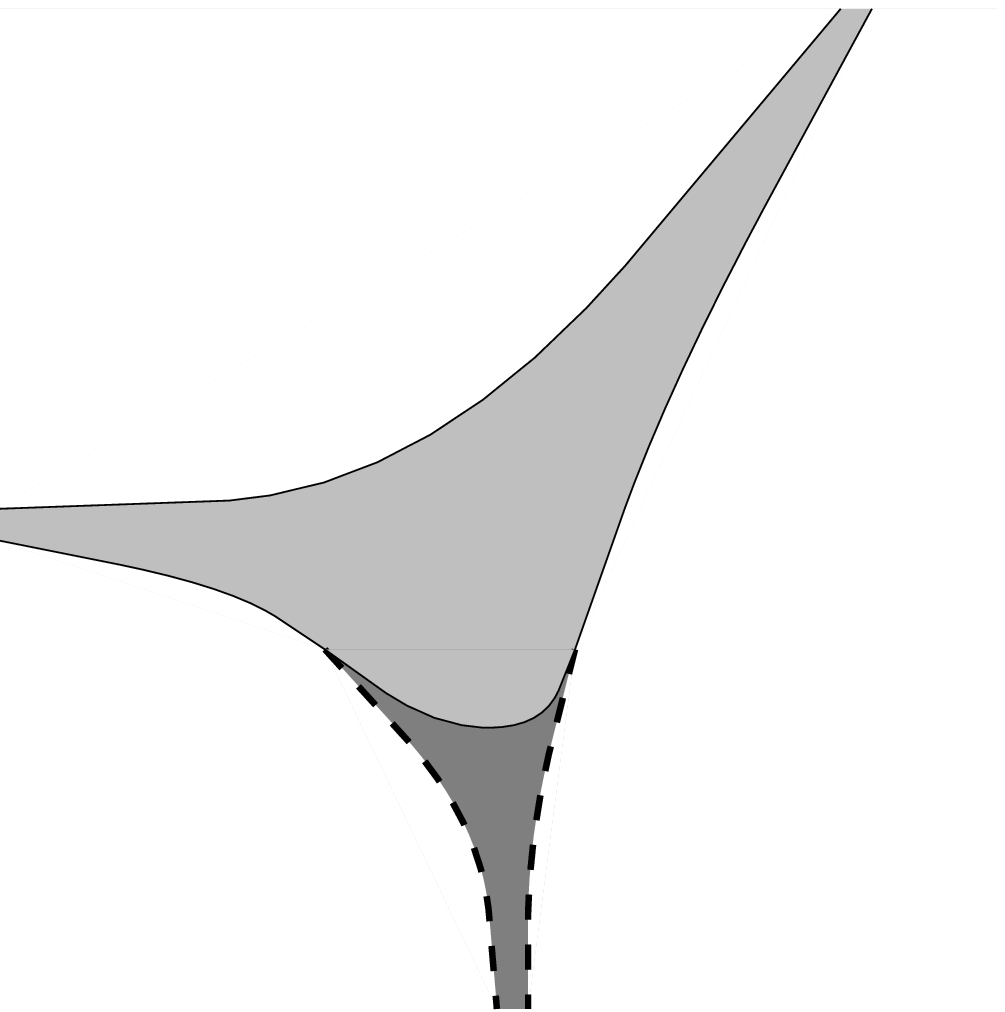 scaled 350}}
\medbreak
\capt{11}{A parabola and its amoeba.}
\enddemo

\vglue-24truept
\section{Proof of the main theorems}
\demo{{\rm 5.1.}\enspace Proof of Theorems {\rm 1} and {\rm 3}}
Theorem 1 is a special case of Theorem~3, where $T=\C P^2$.
All the lemmas in this subsection are formulated under
the hypothesis of Theorem 3.
\enddemo

\specialnumber{5} \proclaim{Lemma}
\label{lemF}
$F=\R A${\rm .}
\endproclaim

\demo{Proof}
By Pick's formula, twice the area of $\Delta$ is equal to twice the
number of lattice points $g$ inside $\Delta$ plus
the number of lattice points $h$ on the boundary minus 2.
Recall that $g$ is the genus of $\C A$ (see \cite{Kh}).
By the maximality assumption (Definition~2)
$g$ is the number of closed components (ovals) of $\R A$.
Note that for any oval $C$ the inverse image ${(\gamma|_C)}^{-1}(x)$,
$x\in\R P^1$, consists at least of two points. 
Indeed, $\mu(C)\subset\R^2$ is also an oval and its projection to $\R$
along the direction determined by $x$ must have two endpoints.

By definition of $d_j$ we have $h=\sum_j d_j$.
By the maximality assumption (Definition~2)
the number of nonclosed components (arcs) of $\R A$
both of whose endpoints belong
to the same axis $l_j$ of $T$ is $\sum_j (d_j-1)=h-n$.
Choose $x\in\R P^1$ to be any point which does not correspond
to the orthogonal direction of the slope
of a side of $\Delta$.
Then for any such arc $B$ the inverse image ${(\gamma|_B)}^{-1}(x)$,
$x\in\R P^1$, consists at least of one point.

Denote by $D$ the $n$ remaining arcs.
The cyclical maximality assumption implies that
${(\gamma|_D)}^{-1}(x)$ consists at least of $n-2$ points
(the sum of the angles of a Euclidean $n$-gon is $\pi(n-2)$).

Adding the above we conclude that ${(\gamma|_D)}^{-1}(x)$
consists at least of $2{\rm Area}\Delta$
points. By Lemma \ref{gauss} it consists precisely 
of $2{\rm Area}\Delta$ points and $\gamma^{-1}(\R P^1)=F=\R A$.
\enddemo

\specialnumber{6} \proclaim{{C}orollary}
\label{noflex}
$\mu(\R A)$ does not have inflection points{\rm .}
\endproclaim

\demo{Proof}
An inflection point of $\mu(\R A)$ would correspond to a real
critical point of $\gamma$.
By Lemma \ref{foldgauss} that would correspond to a singular
point of $F$.
But $\R A$ is nonsingular and thus Lemma \ref{lemF} 
yields a contradiction.
\enddemo

\specialnumber{7} \proclaim{{C}orollary}
The   order of the intersection points $l_j\cap c_j$ on $l_j$
and on $c_j$ {\rm (}\/see Definition~{\rm 2)} agrees{\rm .}
\endproclaim

\demo{Proof}
If these orders do not coincide then one of the arcs of $\mu(\R A)$
must have an inflection point since the ends of $\mu(\R A)$ are convex
in the half-planes cut by the asymptotes.
\enddemo

\specialnumber{8} \proclaim{Lemma}
\label{lemmu}
$\dd\mu(A)=\mu(\R A)${\rm .}
\endproclaim

\demo{Proof}
Lemma \ref{lemF} implies that $\dd\mu(A)\subset\mu(\R A)$.
Suppose that $\dd\mu(A)\neq\mu(\R A)$. Then for some $x\in\R^2$,
$x\notin F$, $\mu^{-1}(x)$ consists of more than 2 points.
Let $L$ be a line passing through $x$ with a rational slope
which is not orthogonal to the slope of a side of $\Delta$.
Let $y\in L$ be a point close to   infinity in $L$ so that
$\mu^{-1}(y)=\emptyset$.
By Corollary \ref{noflex} each component of $\R A$ cuts $\R^2$
into a convex and a nonconvex half.
We call the convex half {\it the interior} of the component
(even if this component is noncompact).
By Lemma \ref{locconv}, if $x$ belongs to the interior of
$a$ components and $y$ belongs to the interior of $b$ components
then the number of points in $\mu^{-1}(x)$ is $2(b-a)$.
But $b=1$ and because of the maximality there is only one arc
which joins the sides of $\Delta$ adjacent to $y$ and only
the interior of this arc may contain $y$. Therefore, $2(b-a)\le 2$. 
\enddemo

\specialnumber{9} \proclaim{{C}orollary}
$\mu|_{\R A}$ is an embedding{\rm .}
\endproclaim

\demo{Proof}
A double point of $\mu(\R A)$ cannot be an inflection
point of a branch by Corollary \ref{noflex}.
Therefore, the other branch
must intersect ${\rm Int}\mu(A)$.
This gives a contradiction to Lemma \ref{lemmu}.
\enddemo

\specialnumber{10} \proclaim{{C}orollary}
For any $(i_1,i_2)\in\Delta$ there exists a unique component
$\Omega_{i_1,i_2}$ of $\R^2-\mu(A)${\rm .} This component
is a disk bounded by an oval of $\mu(\R A)$
if $(i_1,i_2)\in{\rm Int}\Delta$ and a half\/{\rm -}\/plane
bounded by an arc of $\mu(\R A)$ if $(i_1,i_2)\in\dd\Delta${\rm .}
Any component of $\R^2-\mu(A)$ is $\Omega_{i_1,i_2}${\rm ,}
$(i_1,i_2)\in\Delta${\rm .}
\endproclaim

In other words, the amoeba $\mu(A)$ is isotopic to the amoeba
pictured in Figure~12.

\figin{amoeb}{350}
\capt{12}{Amoeba of a curve in cyclically
maximal position in $\R T$.}

\demo{Proof}
The component $\Omega_{i_1,i_2}$ is unique by Theorem 4.
By the same theorem, the value of $\ind$ on any component
of $\R^2-\mu(A)$ belongs to $\Delta$.
But by Lemma \ref{lemmu} any component of $\mu(\R A)$ belongs
to the boundary of a component in $\R^2-\mu(A)$;
because of maximality of $\R\bar{A}$ the number
of components of $\R A$ is equal to the number of lattice points
in $\Delta$.
Thus, we have a disk $\Omega_(i_1,i_2)$
for any $(i_1,i_2)\in{\rm Int}\Delta$ and a half-plane for
any $(i_1,i_2)\in\dd\Delta$.
\enddemo

To reconstruct the topological type of $(\R T,\R\bar{A})$
it now suffices to know the distribution of the components
of $\R A$ among the quadrants of $(\R-0)^2$.
The following lemma determines the lifting of components
of $\mu(\R A)$ to $(\R-0)^2=\mu^{-1}(\R^2)\cap\R T$ and
finishes the proof of Theorem 3 and Theorem 1.

Denote $\R^2_{(-1)^{j_1},(-1)^{j_2}}=
\{(x_1,x_2)\ |\ (-1)^{j_1}x_1>0,\ (-1)^{j_2}x_2>0\}$.
Fix a point $y\in\R A$; then
$\mu(y)\in\dd\Omega_{j_1,j_2}$ for some $(j_1,j_2)\in\Delta$.
Without   loss of generality we may assume that
$y\in\R^2_{(-1)^{j_2},(-1)^{j_1}}$ (otherwise we
change the signs of some of the coordinates).
It turns out that, after this choice of signs,
the same is true for all $\R A\cap\dd\mu(A)$.

\specialnumber{11} \proclaim{Lemma}
\label{keylem}
If $x\in\R A$ and $\mu(x)\in\dd\Omega_{i_1,i_2}$ then
$x\in\R^2_{(-1)^{i_2},(-1)^{i_1}}${\rm .}
\endproclaim

\demo{Proof}
Connect $\mu(x)$ and $\mu(y)$ with a smooth path $Q$ inside $\mu(A)$.
The homology class of the circle $C=(\mu|_A)^{-1}(Q)$ in
$H_1((\C-0)^2)=\Z^2$ is $(j_2-i_2,j_1-i_1)$.
To see this consider two loops $\alpha\subset\mu^{-1}(\mu(x))$ and $\beta\subset\mu^{-1}(\mu(y))$ not intersecting $A$
which represent the same homology class
$(a,b)\in H_1((\C-0)^2)$.
The difference of the linking number of $\alpha$ and $\beta$
with the closure of $A$ in $\C^2$ is equal to
$(j_1-i_1)a-(j_2-i_2)b$.
But on the other hand it must be equal to
the intersection number of an annulus,
connecting $\alpha$ and $\beta$, with $A$.
And that number is equal, in turn,
to the intersection number in $\mu^{-1}(\mu(x))$
of $\alpha$ and the projection of $C$ 
(note that the real 2-torus $\mu^{-1}(\mu(x))$
is a deformational retract of $(\C-0)^2$).

Let $R$ be a path connecting $x$ and $y$ in $\R^2\supset(\R-0)^2$.
Points $x$ and $y$ separate $C$ into two half-circles.
Let $D$ be a membrane in $\C^2$ spanned by the union of $R$
and a half-circle of $C$.
Then $D\cup\conj(D)$ is a membrane for $C$.
The parity of the intersection number
of $D$ with the $x_k$-axis in $\C^2$ ($k=1,2$)
coincides with the parity of the intersection number
of $R$ with the real part of the $x_k$-axis
since all imaginary intersection points come in pairs.
Thus, the intersection number is odd if and only if
$x$ and $y$ are separated by the $x_k$ axis in $\R^2$.
But $[C]=(j_2-i_2,j_1-i_1)\in H_1((\C-0)^2)$ and thus
the linking number of $C$ and the $x_k$-axis is $j_{3-k}-i_{3-k}$.
\enddemo

\demo{{\rm 5.2.}\enspace Proof of Theorem {\rm 2}}
By Remark~3 we may assume that $d\ge 4$.
We deduce Theorem 2 from Theorem 1. Indeed, if $\R\bar{A}$ is in maximal
position with respect to $l_1,l_2,l_3,l_4$ then it is in
maximal position with respect to $l_1,l_2,l_3$ and thus of
the type of Theorem 1. Suppose without   loss of generality
that the $d$ intersection points with $l_4$ belong to the arc $C$
of $\R\bar{A}-(l_1\cup l_2\cup l_3)$ connecting $l_1$ and $l_2$.
Since $\R\bar{A}$ is in maximal position with respect to $l_1,l_2,l_4$
and $d\ge 4$ the region between $C$ and the corner between $l_1$ and $l_2$
must contain an oval of $\R\bar{A}$ by Theorem 1.
But that is impossible by Theorem 1 applied to $\R\bar{A}$
and $l_1,l_2,l_3$.
\enddemo

\bigbreak
\centerline{\bf Appendix. Viro's patchworking and amoebas}
\bigbreak

Viro introduced a patchworking technique for constructing
algebraic curves in his dissertation. This construction is
described in a very elementary way in \cite{IV}.
It turns out that the same construction also provides patchworking for amoebas.

Let us recall the construction. Let $\Delta$ be a convex lattice polygon in $\R^2$,
$T$ be the toric surface corresponding to $\Delta$ and $\mu_T:T\to\Delta$ be
the moment map.
Any function $h:\Delta\cap\Z^2\to\Z$ determines
a subdivision of $\Delta$ into a union of smaller
polygons $\Delta_j$ by the following rule.
Let $R$ be the convex hull of $\{(x,y,t)\in\Z^3\ |\
(x,y)\in\Delta,\ t\ge h(x,y)\}$. The region $R$ is
a semi-infinite polyhedral region. We define $\Delta_j$
to be the vertical projections of its faces.

Let $\{p_j=\sum\limits_{\Delta_j} a_{m,n}x^my^n\}$
be a collection of real polynomials such that
the Newton polygon of $p_j$ is $\Delta_j$.
We assume that the truncations of $p_j$ and $p_k$ to
any common face of $\Delta_j$ and $\Delta_k$
coincide.
Let $\Gamma$ be a side of $\Delta_j$ for some~$j$.
The truncation $p_{\Gamma}$ is a weighted homogeneous
polynomial in two variables.
Suppose that it does not have multiple roots for any $\Gamma$
and all its roots are real.

Let $T_j$ be the toric surface associated to $\Delta_j$
and let $\mu_j:T_j\to\Delta_j$ be the moment map.
Suppose that $p_j$ defines a nonsingular curve $A_j$ in $T_j$.
Suppose that the only singularities of $\mu_j|_{A_j}$
are real folds, imaginary double folds, double cusps
and the two singularities from Remark~9.
Note that any cusp point of $\mu_j|_{A_j}$ is imaginary and
thus, together with the conjugate point, it necessarily
forms a double cusp.
Suppose that the images of the folds
$\mu_j(F_j)\subset{\rm Int}\Delta_j$ of $\mu_j$ do  not
have self-tangencies or triple points.

Define the {\it patchworking} polynomial $p_t$ 
of $\{p_j\}$ by 
$$\sum\limits_{\Delta}a_{m,n}t^{h(m,n)}x^my^n$$
for $t>0$.
The Newton polygon of $p_t$ is $\Delta$.
Denote the zero set of $p_t$ in $(\C-0)^2$ with $A$ and its
zero set in the toric surface $T$ corresponding to $\Delta$
by $\bar{A}$.
Let $F$ be the folds of the moment map $\mu|_A$.

The curve $\bigcup\limits_j\mu_j(F_j)$ is not properly embedded
in ${\rm Int}\Delta$.
To make it proper we need to connect the loose ends on
both sides of $\Gamma_k$ for each interior edge $\Gamma_k$.
We have two pairs of those ends for each zero of
the truncation $p|_{\Gamma_k}$, one pair from each side
of $\Gamma_k$.
The two ends of each pair correspond to a pair of different quadrants
of $(\R-0)^2$. The two ends from the other pair correspond
to the same pair of quadrants; thus we have a natural identification
between the two pairs of ends.
We connect the corresponding ends across $\Gamma_k$
without introducing
a self-intersection or introducing a new ordinary double point
(depending on the mutual position of the two pairs).
Let $f$ be the resulting proper curve in ${\rm Int}\Delta$.
Its isotopy type is determined by the isotopy types
of $\mu_j(F_j)$ and the distribution of their real ends among
the components of $(\R-0)^2$.

\specialnumber{A1}
\proclaim{Proposition}
The pairs $({\rm Int}\Delta,\mu_T(F))$ and 
$({\rm Int}\Delta,f)$ are homeomorphic{\rm .}
Under this homeomorphism the singularities of $\mu_j|_{A_j}$
correspond to the singularities of the same type of $\mu_T|_{A}$
and the new double points of $f$ correspond to the pinching
points of $\mu_T(A)${\rm .}
\endproclaim

\demo{Proof}
Since $R$ is convex there exists a coordinate change
of the type $X=xt^a$, $Y=yt^b$, $a,b\in\Z$ for every
$\Delta_j$ such that the pull-back $P$ of $p$ under
this coordinate change is given by
\begin{eqnarray}
\label{P}
P(X,Y)&=&\sum\limits_{\Delta} a_{m,n}t^{H(m,n)}X^mY^n\\
&=&
t^H(\sum\limits_{\Delta_j}a_{m,n}X^mY^n+
\sum\limits_{\Delta-\Delta_j}a_{m,n}t^{H(m,n)-H}X^mY^n),
\nonumber
\end{eqnarray}
$H(m,n)\ge H$ and $H(m,n)=H$ if and only if $(m,n)\in\Delta_j$.
On the other hand, the pull-back of $p_j$ is
\begin{equation}
\label{Pj}
P_j(X,Y)=\sum\limits_{\Delta_j}a_{m,n}t^{H(m,n)}X^mY^n.
\end{equation}
Let $D_j$ be a bounded disk region in the $(X,Y)$-plane
such that the amoeba of $P_j$ is essentially contained in it;
i.e., $\mu_j^{-1}(\R^2-D_j)\cap\{(X,Y)\ |\ P(X,Y)=0\}$ is
a disjoint union of annuli which project under $\mu$
to strips connecting  infinity with the  $D$ where the
folds are the only singularities of the projection.
Let $K_j$ be the region corresponding to $D_j$
in the $(x,y)$-plane. Note that $K_j$ depends on the
value of $t$.

If $t>0$ is sufficiently small then the zero sets of $p$
and $p_j$ in $\mu^{-1}(K_j)\subset(\C-0)^2$ are sufficiently
close (cf. \eqref{P} and \eqref{Pj}).
On the other hand for sufficiently small $t>0$ and $j\neq l$,
$K_j\cap K_l=\emptyset$.
By assumption the singularities of $\mu|_{A_j}$
are from our list; they are stable under real deformation.
Therefore, $\mu|_A$
has the same singularities over $K_j$ as $\mu|_{A_j}$.
By the other assumption $\mu(F_j)$ does not have double tangencies
or triple points. Therefore, $\mu(F)\cap K_j$ and $\mu(F_j)\cap K_j$
are isotopic in $K_j$.

Note that we chose $K_j$ so that $A_j-\mu^{-1}(K_j)$ is
a collection of disjoint annuli.
Recall that the genus of $A$ is equal to the number of the lattice
points in ${\rm Int}\Delta$ and the same is true for
$A_j$. This implies that $A-\bigcup\limits_j\mu^{-1}(K_j)$ is
a collection of disjoint annuli.
By Lemma \ref{foldgauss} the image of each of these annuli
double covers a small disk around
the orthogonal direction to $\Gamma_j$ in $\C P^1$
under the logarithmic Gauss map.
Therefore, this covering has two branching points.
If they are real we have a pinching. If they are imaginary
then the image of the folds of the annuli under $\mu$ is
embedded.
\enddemo

The curve obtained as a result of patchworking
of polynomials whose Newton polygons $\Delta_j$
are triangles of area $\frac{1}{2}$ in $\R^2$
are called T-curves (see \cite{IV}).
Note that any strictly convex function
$h:\R^2\to\R$ gives a decomposition of $\Delta$ into
such triangles.

\specialnumber{A2}
\proclaim{{C}orollary}
If $A$ is a {\rm T-}\/curve then for every $(m,n)\in\Delta\cap\Z^2$
there exists a component $\Omega$ of $\R^2-\mu(A)$ such that
$\ind\Omega=(m,n).$
\endproclaim

Indeed, Proposition A1 implies this even in a slightly more
general situation when ${\rm Ind}\Delta_j\cap\Z^2=\emptyset$
for each $j$.

Let $\Delta_1$ be the convex hull of $(m_0,n_0)$,
$(m_1,n_1)$ and $(m_2,n_2)$.
Let $\Delta_2$ be the convex hull of $(m_1,n_1)$,
$(m_2,n_2)$ and $(m_3,n_3)$.
Suppose that ${\rm Area}(\Delta_1)=
{\rm Area}(\Delta_2)=\frac{1}{2}$ and
${\rm Int}\Delta_1\cap{\rm Int}\Delta_2=\emptyset$.
Suppose that $\Delta_1$ and $\Delta_2$ are contained in
a patchworking of a Newton polygon $\Delta$ and
$a_j\neq 0$ are the coefficients at $x^{m_j}y^{n_j}$.
Let $\Gamma=\Delta_1\cap\Delta_2$.

\specialnumber{A3}
\proclaim{Lemma}
If $a_1a_2<0$ then
pinching at $\Gamma$ occurs if and only if
$a_0a_3<0${\rm .}
If $a_1a_2>0$ and $(m_0,n_0)\equiv(m_2,n_2)\!\!\pmod{2}$ then
pinching at $\Gamma$ occurs if and only if
$a_0a_3<0${\rm .}
If $a_1a_2>0$ and $(m_0,n_0)\not\equiv(m_2,n_2)\!\!\pmod{2}$ then
pinching at $\Gamma$ occurs if and only if
$a_0a_3>0${\rm .}
\endproclaim

Note that this lemma agrees with Example 1.
\demo{Proof}
The pull-back of $a_jx^{m_j}y^{n_j}$ under the sign change
of the $x$-coordinate is $(-1)^{m_j}x^{m_j}y^{n_j}$.
Such a change does not affect amoeba.
Thus, changing signs of $x$ or $y$ if needed we
may assume that $a_1a_2<0$.
In this case the curve
$a_0x^{m_0}y^{n_0}+a_1x^{m_1}y^{n_1}+a_2x^{m_2}y^{n_2}=0$
(which is just a reparametrization of a line in $(\R-0)^2$,
since ${\rm Area}(\Delta_1)=\frac{1}{2}$)
intersects the positive quadrant in $(\R-0)^2$.
The image of this intersection under $\mu$ does not
have inflection points. Its curvature is directed towards
$(m_1,n_1)$
if $a_0a_2>0$ and towards
$(m_2,n_2)$ otherwise. The same is true for
$a_3x^{m_3}y^{n_3}+a_1x^{m_1}y^{n_1}+a_2x^{m_2}y^{n_2}=0$
in $\Delta_2$. Thus pinching occurs if and only if $a_0a_3<0$.
\enddemo

\specialnumber{A4}
\proclaim{{C}orollary}
For any convex lattice polygon $\Delta$
there exists a polynomial $p$ whose
Newton polygon is $\Delta$ and
the corresponding real curve $\R\bar{A}\subset\R T$
is in cyclically maximal position{\rm .}
\endproclaim

\demo{Proof}
We construct $\R\bar{A}$ as a T-curve.
Decompose $\Delta$ into a union of triangles
$\Delta_j$ of area $\frac{1}{2}$ using a strictly
convex function $h$.
Take $a_{m,n}=(-1)^{(m-1)(n-1)}$ (cf. \cite{IV})
as the coefficients of $p_j$.
By Lemma A3, pinching does not occur. Thus, by Proposition A1,
the corresponding curve is in cyclically maximal 
position in $\R T$.
\enddemo

\bye